\documentclass{article}

\usepackage{arxiv}

\usepackage{dsfont}
\usepackage{amssymb}
\usepackage{amsmath}
\usepackage{amsfonts}
\usepackage{mathtools}
\usepackage{empheq}
\usepackage[utf8]{inputenc} 
\usepackage[T1]{fontenc}    
\usepackage{hyperref}       
\usepackage{url}            
\usepackage{booktabs}       
\usepackage{amsfonts}       
\usepackage{nicefrac}       
\usepackage{microtype}      
\usepackage{cleveref}       
\usepackage{lipsum}         
\usepackage{graphicx}
\usepackage{doi}

\newtheorem{lemma}{Lemma}

\newtheorem{theorem}{Theorem}
\newtheorem{corollary}{Corollary}
\newtheorem{remark}{Remark}

\title{Existence of an attractor and Horseshoe in multidimensional H\'{e}non map}

\author{D.~A.~Grechko \\
	Volga State University of Water Transport\\
	Nizhny Novgorod, 603950 Russia \\
        Lobachevsky State University of Nizhny Novgorod \\
        Nizhny Novgorod, 603022 Russia \\
	\texttt{d.grechko.18@gmail.com} \\
 \And
 V.~N.~Belykh \\
	Volga State University of Water Transport\\
	Nizhny Novgorod, 603950 Russia \\
        Lobachevsky State University of Nizhny Novgorod \\
        Nizhny Novgorod, 603022 Russia \\
	\texttt{Belykh.vn@vsuwt.ru} \\
	\And
	N.~V.~Barabash \\
	Volga State University of Water Transport\\
	Nizhny Novgorod, 603950 Russia \\
        Lobachevsky State University of Nizhny Novgorod \\
        Nizhny Novgorod, 603022 Russia \\
	\texttt{barabash@itmm.unn.ru} \\
}

\renewcommand{\shorttitle}{Existence of an attractor and Horseshoe in multidimensional H\'{e}non map}

\hypersetup{
pdftitle={Existence of an attractor and Horseshoe in multidimensional generalized H\'{e}non map},
pdfsubject={math.DS},
pdfauthor={V.~N.~Belykh,N.~V.~Barabash,D.~A.~Grechko},
pdfkeywords={dynamical systems, nonlinear maps, attractors, hyperbolicity, bifurcations},
}
\begin{document}
\maketitle

\begin{abstract}
In this paper using approach of 1-D auxiliary maps we prove the existence of trapping domains containing attractors of the multidimensional H\'{e}non-like maps.
For both of quadratic and cubic nonlinearities we obtain sufficient conditions of topological Smale horseshoes existence.
The complex structure of attractors is discussed in the case of small coupling parameter.
\end{abstract}

\keywords{dynamical systems \and nonlinear maps \and attractors \and H\'{e}non map \and bifurcations}

\section{Introduction}

Strange attractors as attracting invariant sets of entire unstable trajectories can be divided into 3 types: hyperbolic, whose structure does not change at all points of the interval of the parameter characterizing the deformations of the dynamical system; singular hyperbolic structures whose structure changes only at bifurcation points; quasistranger - strange not on an interval, but on a point, usually Cantor set of parameters.

In the class of discrete-time dynamical systems (mappings), examples of hyperbolic attractors are Anosov diffeomorphisms \cite{Anosov1963}, examples of singular-hyperbolic attractors are the Lorentz attractor for model mappings \cite{AfraimovichShilnikov1982,Robinson1989}, Lozi \cite{Lozi1978} and Belykh \cite{BelykhMatSbornik1995}, and others. The conditions for the existence of reduced attractors are written out analytically by virtue of a specific form of mappings.

We consider a $(n+1)$-dimensional map of the form
\begin{equation}
T:\quad \begin{array}{lcl}
\bar{x}=&f(x)+\mathds{1}y&\triangleq \; g(x,y),\\
\bar{y}=&\textbf{\textrm{b}}x+Ay&\triangleq \; L(x,y),
\end{array}
\label{map:T}
\end{equation}
where $x\in\mathbb{R}^1$, $y=\textrm{column}(y_1,y_2,\ldots, y_n)\in \mathbb{R}^n$, $\mathds{1}$ is the all ones line of length $n$, $f(x)$ is a continuous smooth or piecewice function, $\textbf{\textrm{b}}=\textrm{column}(b,0,0,\ldots,0)$,  $b$ is the parameter, $A=[a_{ij}]_n^n$ is the normalized lower shift $(n\times n)$-matrix with entries $a_{ij}=a_j\delta_{i,j+1}$, where $a_j$ are parameters and $\delta_{i,j+1}$ is the Kronecker delta symbol. 
The bar symbol over variables as usually denotes the next iteration.  

The change of the variables $y_i=a_i v_i$, $i=1,\ldots,n$, and the parameters $b=q_1$, $a_{i}=\dfrac{q_{i+1}}{q_{i}}$, $i=1,\ldots,n-1$, transforms the map (\ref{map:T}) to the form of generalized H\'{e}non-like map \cite{gonchenko2005three,Li_2006}  

\begin{equation*}
\begin{array}{lcl}
\bar{x}=f(x)+\sum\limits_{j=1}^n q_j v_j,\\
\bar{v}_1=x, \\
\bar{v}_{i+1}=v_{i},\quad i=1,\ldots, n-1.
\end{array}
\label{sys:Henon_1}
\end{equation*} 

Note, that in two-dimensional case for
\begin{equation}
n=1 \textrm{ and } f(x)=\mu-x^2
\end{equation}
this map becomes original H\'{e}non map \cite{henon1976two}.

Michel H\'{e}non proposed this map \cite{henon1976two} as an abstract example of a dynamical system with a strange attractor. Nevertheless, it can serve to describe the dynamics of a number of physical systems, for example, a dissipative oscillator under the influence of an external force, the magnitude of which depends nonlinearly (quadratically) on the coordinate of the oscillator \cite{Franks1969AnosovDiff} and the rotator under a pulsed periodic action.

The original 2-D H\'{e}non map and its multidimensional analog was considered in set of papers \cite{CarlesonBenedicks1991,GonchenkoMeissOvsyannikov2006,Heagy1992,SterlingDullin1999,GonchenkoPhysica2016}. The authors studied the main hard problem whether (when) H\'{e}non attractor is strange or not.

New type of attractor of multidimensional H\'{e}non map has been discovered in the papers \cite{GonchenkoLiMalkin2008}. This attractor is similar in appearance to the Lorenz attractor thought its structure is different. In \cite{kennedy2001topological} the existence of Smale horseshue was proved and ergodic properties were studied in the case of multidimensional H\'{e}non map closed to 1-D unimodal map.

In this paper we consider localization of attractor and existence of Smale horseshue in multidimensional H\'{e}non map without small parameters.

We assume that the parameters of the map ($\ref{map:T}$) satisfy the conditions
\begin{equation}
|b|<1,\quad |a_i|<a<1,\quad i=1,\ldots,n-1.
\label{eq:a<1_b<1} 
\end{equation}

We use the Manhattan norm 
\begin{equation}
\Vert y \Vert=\sum\limits_{i=1}^n |y_i|.
\label{eq:norm}
\end{equation}

\section{Auxiliary 1-D map}

We introduce two domains $D_{\alpha}=D_x\times D_y$ and 
$D=\bar{D}_x\times D_y$, where 
$D_x=\lbrace \alpha^-<x<\alpha^+ \rbrace$, 
$\bar{D}_x=\lbrace x\in\mathbb{R}^1 \rbrace$ and 
$D_y=\lbrace \Vert y\Vert<\gamma \rbrace$, for some constants $\alpha^-, \alpha^+, \gamma>0$.
Denote $\alpha=\max\lbrace |\alpha^-|,|\alpha^+|\rbrace$.

\begin{lemma}
Let parameter $\gamma$ defining the domain $D_y$ be
\begin{equation}
\gamma=\dfrac{\sigma|b|}{1-a}.
\label{eq:gamma}
\end{equation}

Then image $TD_{\alpha}$ lies in $D$, i.e. 
\begin{equation*}
    TD_{\alpha}\subset D.
\end{equation*}
\label{lem:lemma_1}
\end{lemma}

\textbf{Proof.} 
According to (\ref{map:T}), (\ref{eq:norm}) the norm of the $y$-vector satisfies the next inequality
\begin{equation}
\Vert \bar{y} \Vert \leq 
|b||x| + \sum\limits_{i=1}^{n-1} |a_{i}||y_{i}|.
\label{eq:y_norm}
\end{equation}

Due to condition (\ref{eq:a<1_b<1}) for any coordinate $x$ of preimage lying in $D_x$,  $x\in D_x$ using (\ref{eq:y_norm}) we obtain inequalities 
\begin{equation}
\Vert \bar{y} \Vert <\alpha |b| + a(\Vert y \Vert - |y_n|)<\alpha |b|+a\Vert y \Vert.
\label{eq:norm_ineq}
\end{equation}

From this condition it follows that for any point $(x,y)$ from domain $\lbrace \Vert y \Vert>\dfrac{\alpha |b|}{1-a}=\gamma, x\in D_x \rbrace$ inequality $\Vert \bar{y} \Vert < \Vert y \Vert$ is valid.
This implies that $TD_\alpha\subset D$ for any $(x,y)\in D_\alpha$. 
$\Box$

\begin{corollary}
From Lemma~\ref{lem:lemma_1} it follows, that for
$(x,y)\in D_{\alpha}$, i.e. for $\alpha^-<x<\alpha^+$,
$\Vert y \Vert<\gamma$, image $(\bar{x},\bar{y})=T(x,y)$ satisfies conditions
\begin{equation}
    f(x)-\Vert y \Vert<\bar{x}<f(x)+\Vert y \Vert,\quad \Vert y \Vert<\gamma.
    \label{eq:corol_T(xy)}
\end{equation}
\end{corollary}

Introduce two auxiliary 1-D maps
\begin{equation}
\begin{array}{ll}
x\mapsto f(x)+\gamma\triangleq  f^+(x),\\
\\
x\mapsto f(x)-\gamma\triangleq  f^-(x),
\end{array}
\label{sys:comp_x}
\end{equation}
playing an essential role in the study of map (\ref{map:T}).

\begin{corollary}
For any $(x,y)\in D_{\alpha}$ the image
coordinates satisfy conditions
\begin{equation}
\begin{array}{c}
f(x)-\gamma=f^-(x) < \bar{x} < f^+(x)=f(x)+\gamma,\\
\\
     \bar{y}\in D_y.
\end{array}
\label{eq:statement_2_cond}
\end{equation}
\label{statement_aux_sys}
\end{corollary}

\bigskip

\section{Existence of attractor}

As far as $D_{\sigma}\subset D$ from Lemma~\ref{lem:lemma_1} it follows that for existence of map~\eqref{map:T} attractor it is sufficient ti find such conditions on function $f(x)$ and interval $D_x$ such that $TD_{\sigma}\subset D_{\sigma}$.

\begin{theorem}
    
    Let function $f(x)$ satisfy the next condition
    \begin{equation}
    \begin{array}{c}
        \alpha^-+\gamma\leq f(x)\leq \alpha^+-\gamma,\\
         \\
         x\in D_x,
    \end{array}
    \label{eq:theorem_1_cond}
    \end{equation}
the the map~\eqref{map:T} has an attractor $A\subset D_{\alpha}$.
\label{th:attractor_gamma}
\end{theorem}

\textbf{Proof.}
Let $(x,y)$ be any point lying in domain $D_{\alpha}$.
Then coordinates of image $(\bar{x}, \bar{y})=T(x,y)$ due to Corollary~\ref{statement_aux_sys} lie in $D_{\sigma}$ if $f(x)-\gamma>\alpha^-$ and $f(x)+\gamma<\alpha^+$ what is equivalent to inequalities~\eqref{eq:theorem_1_cond}.
Hence, $TD_{\sigma}\subset D_{\sigma}$ and map~\eqref{map:T} has an attractor $A=\lim\limits_{n\rightarrow\infty} T^nD$. $\Box$

\begin{remark}
    For fixed values $\alpha^-$ and $\alpha^+$, Theorem~\ref{th:attractor_gamma} condition \eqref{eq:theorem_1_cond} shows the class of functions $f(x)$ for which map~\eqref{map:T} has an attractor.
    For given function $f(x)$ satisfying \eqref{eq:theorem_1_cond} values $\alpha^-$ and $\alpha^+$ are undefined and the problem of attractor finding via Theorem~\ref{th:attractor_gamma} leads to the problem of these values existence. 
\end{remark}

Consider this problem. 
Note that map~\eqref{map:T} for $b=0$ has 1-D manifold $\lbrace y_1=y_2=\cdots =y_n=0\rbrace$ dynamics on which is defined by 1-D map $x\mapsto f(x)$.
Let function $f(x)$ have invariant interval $I=\lbrace \alpha_0^-<x< \alpha_0^+ \rbrace$, $fI\subset I$.
Then attractor of map~\eqref{map:T} for $b=0$ lies at interval $I$ and satisfies condition (*) for $\gamma=0$, $\alpha^{\pm}=\alpha_0^{\pm}$.
Denote extreme values of function $f(x)$ at interval $I$ 
\begin{equation}
\begin{array}{lcl}
x^-=\inf\limits_{x\in I} f(x),\\
\\
x^+=\sup\limits_{x\in I} f(x).
\end{array}
\label{eq:x+_x-}
\end{equation}

\begin{lemma}
Let for given $\gamma$ extreme values \eqref{eq:x+_x-} of function $f(x)$ satisfy inequalities
\begin{equation}
    \begin{array}{c}
         x^--\gamma>\alpha_0^-, \\
         \\
         x^++\gamma<\alpha_0^+.
    \end{array}
    \label{eq:lemma_attr_D_a}
\end{equation}
Then map~\eqref{map:T} for $b>0$ has an attractor lying in domain $D_{\alpha}$ with $\alpha^-=\alpha_0^-$ and $\alpha^+=\alpha_0^+$.
\label{lem:D_a_invariant}
\end{lemma}

\textbf{Proof.} 
Putting $\alpha^{\pm}=\alpha_0^{\pm}$ we choose $D_x=I$.
Then at \eqref{eq:lemma_attr_D_a} condition of Theorem~\ref{th:attractor_gamma} holds and map~\eqref{map:T} has an attractor. $\Box$ 

\section{Attractors' existence of H\'{e}non map}

\subsection{The case of quadratic nonlinearity}

Consider the case of quadratic nonlinearity  
\begin{equation}
f(x)=\mu-x^2,
\label{eq:parabola}
\end{equation}
which corresponds to the original 2-D H\'{e}non map.
We introduce the detailed proof of the next
\begin{theorem}
Let the following condition hold
\begin{equation}
0<\mu<2\left(1-\dfrac{|b|}{1-a}\right)^2.
\label{eq:statement_1_condition_mu<2}
\end{equation}

Then the multidimensional H\'{e}non map~\eqref{map:T}, \eqref{eq:parabola} has invariant domain 
\begin{equation}
D=\left\lbrace |x|<\alpha=\dfrac{\mu(1-a)}{1-a-|b|}, \Vert y\Vert <\gamma=\dfrac{\mu|b|}{1-a-|b|}\right\rbrace,
\label{eq:statement_1_D}
\end{equation}
with an attractor $A$.
\label{theorem_logistic_Henon_attractor}
\end{theorem}

\begin{figure}
\centering
\includegraphics[width=0.4\linewidth]{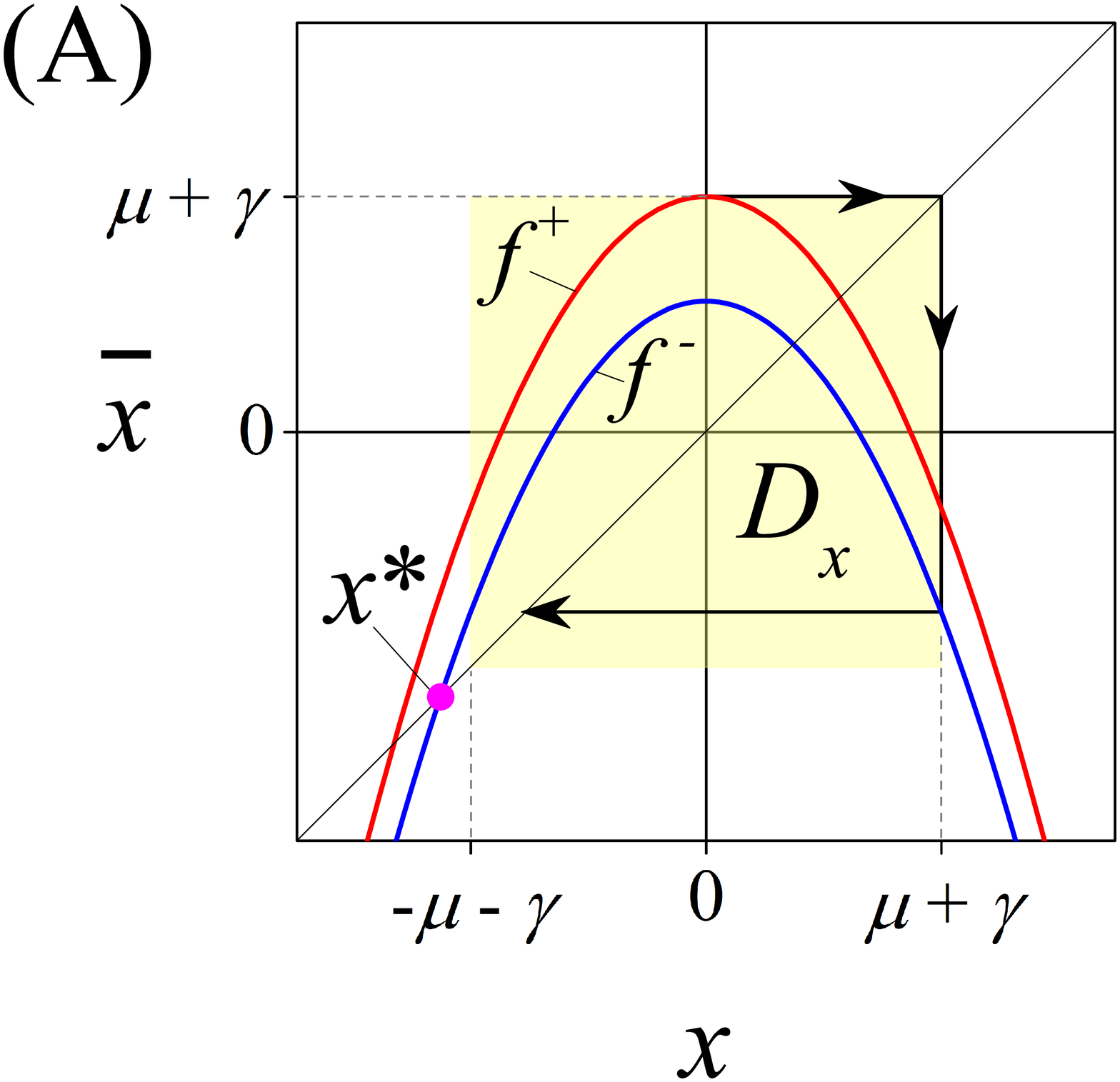}
\includegraphics[width=0.4\linewidth]{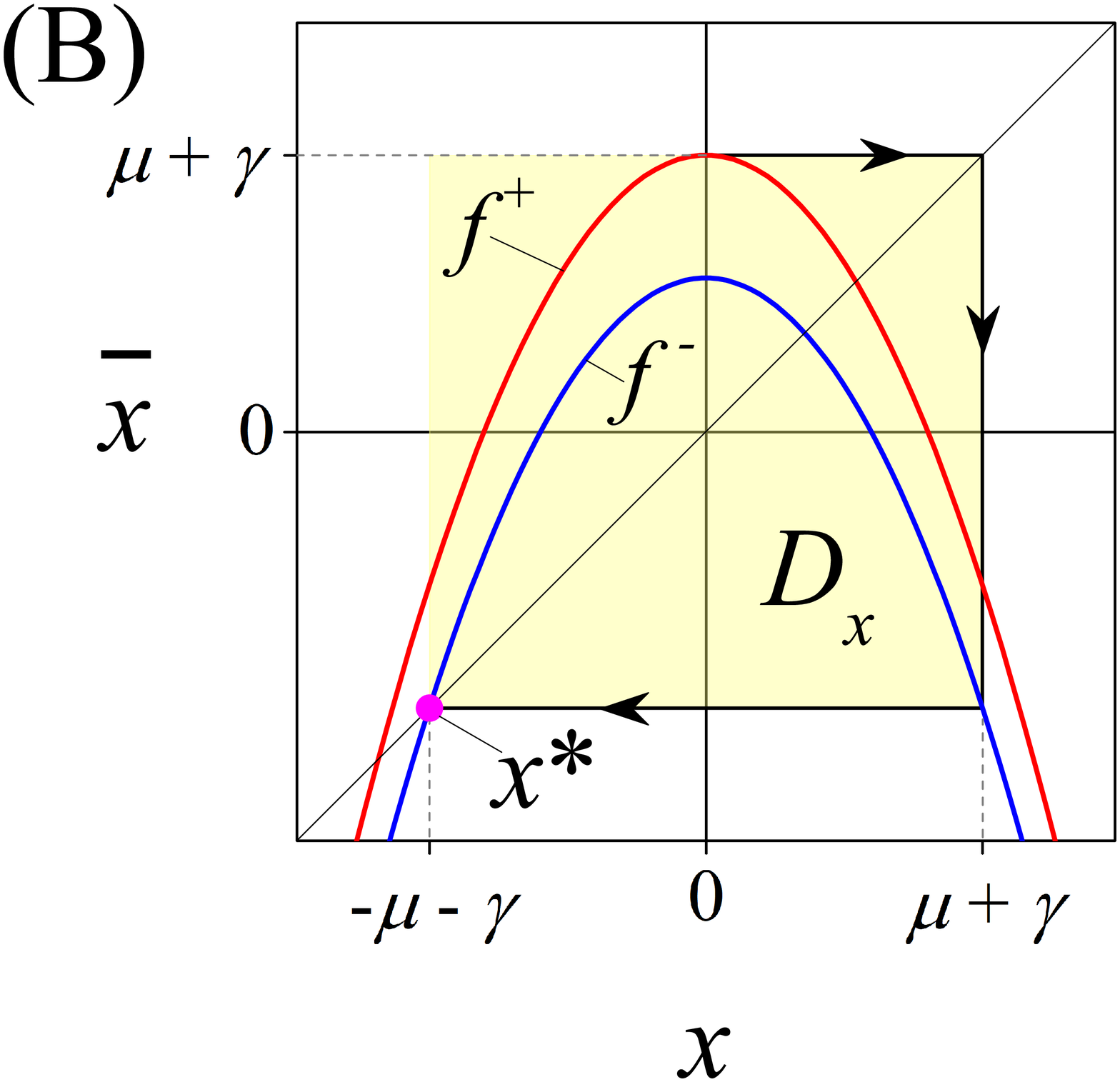}
\caption{Illustration of Theorem's~\ref{theorem_logistic_Henon_attractor} conditions for parameters (A) $\mu=1.028$, $\gamma=0.294$, $x^*=-1.492$; (B) $\mu=1.21$, $\gamma=0.346$, $x^*=-1.556$.}
\label{fig:logistic_aux}
\end{figure}

\textbf{Proof.} 
The parameter $\gamma$ (\ref{eq:gamma}) due to Lemma~\ref{lem:lemma_1} defining the domain $D_y$ for any function $f(x)$ including (\ref{eq:parabola}) depends on the parameter $\alpha$. 
Now we obtain the parameter $\alpha$ defining the interval $D_x$ in the case of the quadratic function (\ref{eq:parabola}).
For this we rewrite 1-D auxiliary maps (\ref{sys:comp_x}) in the following form (see Fig.~\ref{fig:logistic_aux})
\begin{equation}
\begin{array}{l}
g^+:\quad \bar{x}=\mu-x^2+\gamma,\\
g^-:\quad \bar{x}=\mu-x^2-\gamma.
\end{array}
\label{eq:proof_quadr_aux}
\end{equation}

Due to the assumption of Theorem~\ref{th:attractor_gamma}
the expression for $\alpha^+$ takes the form
\begin{equation}
\alpha^+=\mu+\gamma.
\label{eq:a=mu+gamma}
\end{equation}

We choose $\alpha^-=g^-(\alpha^+)$ as the left boundary of $D_x$ being the minimal value of the function $g^-(x)$ on the interval $(\alpha^-, \alpha^+)$.

The condition $g^+(\alpha^+)<\alpha^+$ holds for any positive parameters $\mu$ and $\gamma$,
and the condition $g^-(\alpha^-)>\alpha^-$ is true for  
\begin{equation}
x_l<\alpha^-<x_r,
\label{eq:proof_x^2_ineq}
\end{equation} 
where $x_{l,r}=\frac{-1\mp\sqrt{1+4(\mu-\gamma)}}{2}$ are left and right fixed points of the map $g^-$.
From this inequality we obtain the condition
\begin{equation}
(\mu+\gamma)^2-2\mu<0,
\label{eq:proof_therem_2_condition}
\end{equation}
providing not only the invariance of the interval $(\alpha^-, \alpha^+)$, 
but also the invariance of the extended interval $D_x=(-\alpha, \alpha)$, where $\alpha=\alpha^+$
because the left inequality in (\ref{eq:proof_x^2_ineq}) can be changed by the inequality $x_l<-\alpha<\alpha^-$, providing the invariance of $(-\alpha, \alpha)$.

Moreover, under the condition (\ref{eq:proof_therem_2_condition}) due to (\ref{eq:gamma}) we obtain
\begin{equation}
\gamma=\dfrac{\alpha|b|}{1-a}=\dfrac{\mu|b|}{1-a-|b|}.
\label{eq:gamma_proof}
\end{equation}

Substituting this parameter $\gamma$ in (\ref{eq:a=mu+gamma}) we obtain the value
\begin{equation}
\alpha=\mu\left(1+\dfrac{|b|}{1-a-|b|}\right),
\end{equation}
depending only on the parameters of map~\eqref{map:T}, \eqref{eq:parabola}.

The condition (\ref{eq:statement_1_condition_mu<2}) follows from inequality (\ref{eq:proof_therem_2_condition}), where $\gamma$ is given by (\ref{eq:gamma_proof}).

Therefore, according to Theorem~\ref{th:attractor_gamma}, the map $T$ with quadratic nonlinearity $f(x)$ (\ref{eq:parabola}) has an attractor $A$, lying in the domain $D$ with boundaries $\alpha$ and $\gamma$ explicitly defined via parameters of map~\eqref{map:T}, \eqref{eq:parabola}. $\Box$

\subsection{The case of cubic nonlinearity}

Consider the case of odd cubic nonlinearity
\begin{equation}
f(x)=x^3-\mu x.
\label{eq:qubic}
\end{equation}

\begin{theorem}
Let the following condition be fulfiled 
\begin{equation}
|b|<\dfrac{(1-a)(3-\mu)}{3}, \quad 0<\mu<3,
\label{eq:theorem_qubic_condition}
\end{equation}
Then map~\eqref{map:T}, \eqref{eq:qubic} has an attractor $A$ in the invariant domain 
\begin{equation}
D=\left\lbrace |x|<\alpha=2\sqrt{\dfrac{\mu}{3}},\quad
\Vert y\Vert <\gamma=2\left(1-\dfrac{\mu}{3}\right)\sqrt{\dfrac{\mu}{3}}\right\rbrace.
\label{eq:theorem_qubic_domain}
\end{equation}
\label{theorem_cubic}
\end{theorem}

\begin{figure}
\centering
\includegraphics[width=0.4\linewidth]{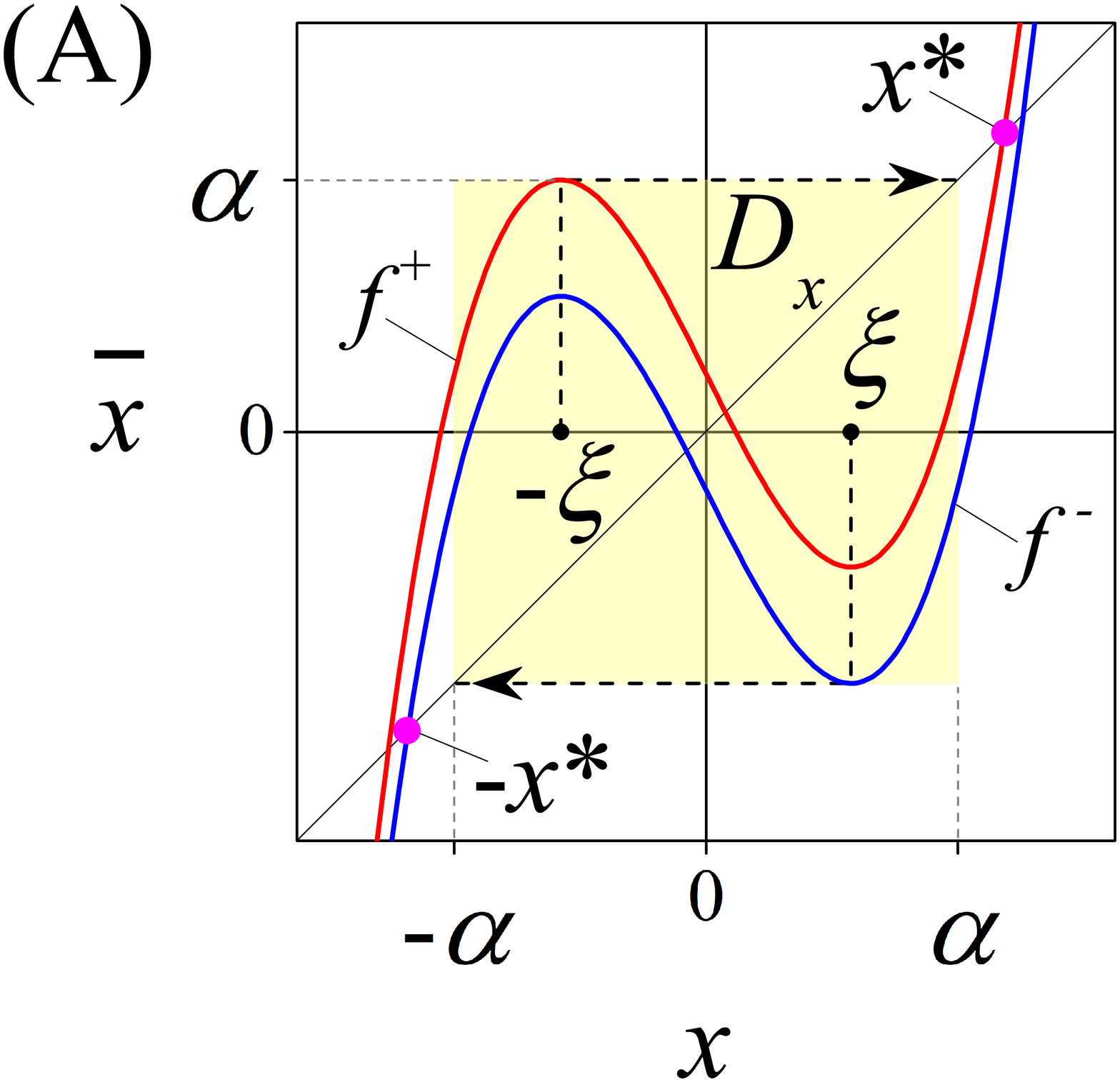}
\includegraphics[width=0.4\linewidth]{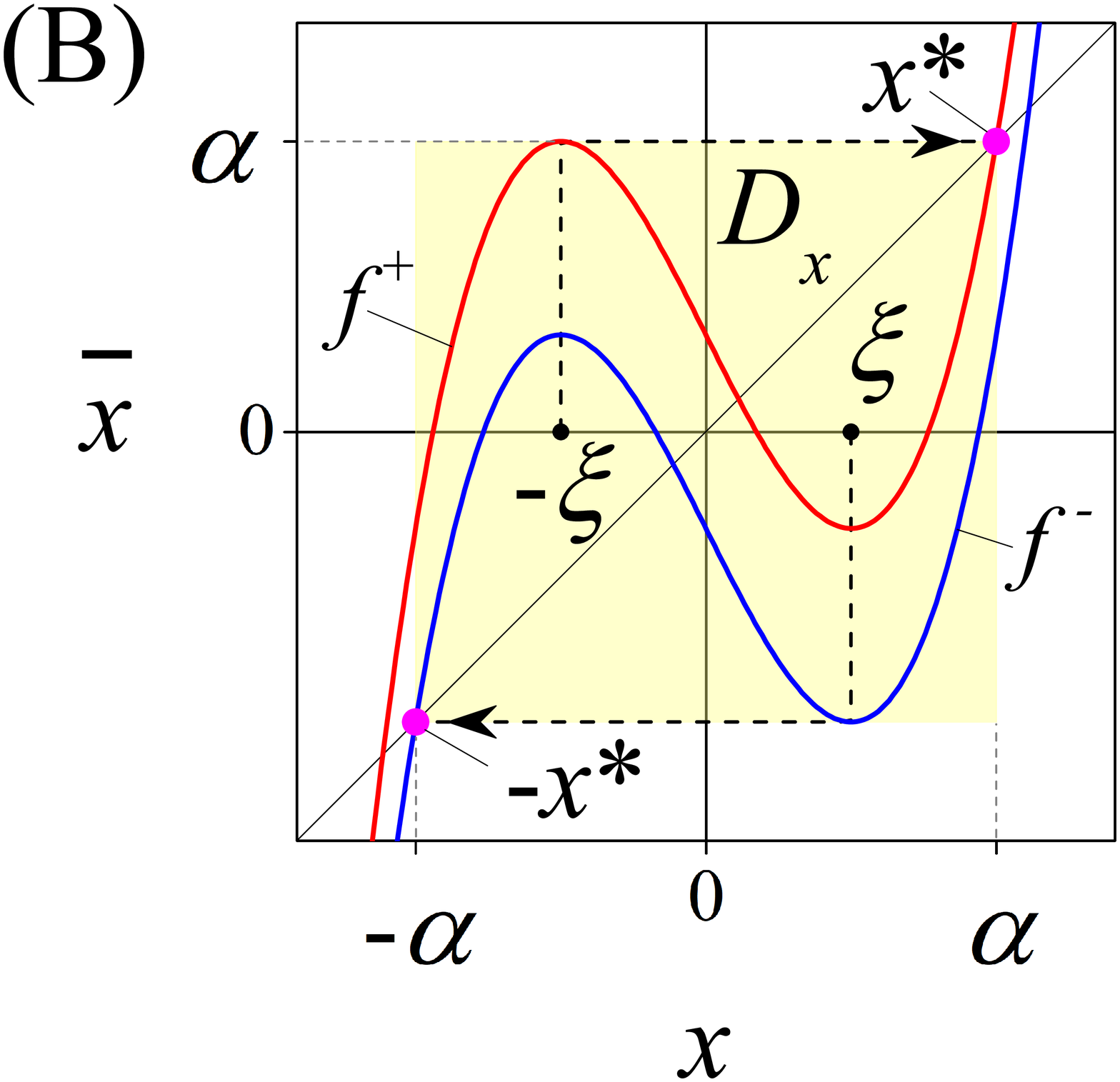}
\caption{Illustration of Theorem's~\ref{eq:theorem_qubic_condition} conditions for parameters (A) $\mu=2$, $\gamma=0.326$, $x^*=1.675$; (B) $\mu=2$, $\gamma=0.544$, $x^*=1.633$.}
\label{fig:cubic_aux}
\end{figure}

\textbf{Proof.} 
Similarly to (\ref{eq:proof_quadr_aux}) we introduce 1-D auxiliary maps in the form (see Fig.~\ref{fig:cubic_aux})
\begin{equation}
\begin{array}{lcl}
g^+:\quad \bar{x}=x^3-\mu x + \gamma,\\
g^-:\quad \bar{x}=x^3-\mu x - \gamma.\\
\end{array}
\label{eq:proof_theorem_quadr_aux}
\end{equation}

The function $f(x)$ \eqref{eq:qubic} has two extreme points, namely the maximum $\bar{x}^+=f(-\xi)$ and the minimum $\bar{x}^-=f(\xi)$, where $\xi=\sqrt{\dfrac{\mu}{3}}$.
Hence according to (\ref{eq:proof_theorem_quadr_aux}) we choose $\alpha^+$ and $\alpha^-$ as
\begin{equation}
\begin{array}{lcl}
\alpha^+=-\xi^3+\mu\xi+\gamma,\\
\alpha^-=\xi^3-\mu\xi-\gamma.
\end{array}
\label{eq:proof_cubic_alpha+-}
\end{equation}

We define the interval $D_x=(-\alpha, \alpha)$ where $\alpha=\alpha^+=|\alpha^-|$.
Substituting the expression (\ref{eq:proof_cubic_alpha+-}) for the parameter $\alpha$ in the expression for the parameter $\gamma$ (\ref{eq:gamma}) we obtain 
\begin{equation}
\gamma=\xi(\mu-\xi^2)\beta=\dfrac{2\xi\mu}{3}\beta,
\label{eq:proof_cubic_gamma}
\end{equation} 
where 
\begin{equation}
\beta=\dfrac{|b|}{1-a-|b|},
\label{eq:proof_beta}
\end{equation}
and hence from (\ref{eq:proof_cubic_alpha+-}) we get
\begin{equation}
\alpha=\dfrac{2\xi\mu}{3}(1+\beta).
\label{eq:proof_cubic_alpha}
\end{equation}

From inequalities (\ref{eq:theorem_qubic_condition}), (\ref{eq:proof_theorem_quadr_aux}) it follows that the condition of the invariance of the interval $D_x$  takes the form
\begin{equation}
\alpha^3-\mu\alpha+\gamma<\alpha.
\label{eq:proof_cubic_condition_1}
\end{equation}

Substituting the expressions  (\ref{eq:proof_cubic_gamma}) and (\ref{eq:proof_cubic_alpha}) for the parameters $\gamma$ and $\alpha$ in the inequality (\ref{eq:proof_cubic_condition_1}) we obtain
the inequality
\begin{equation}
p^3-\dfrac{27}{4}p-\dfrac{27}{4}<0,
\label{eq:proof_theorem_cubic_ineq_p}
\end{equation}
where $p=\mu(1+\beta)>0$.

It turns out that the inequality (\ref{eq:proof_theorem_cubic_ineq_p}) is equivalent to the inequality
\begin{equation}
(p-3)(p+1.5)^2<0.
\end{equation} 
from which it follows that $p=\mu(1+\beta)<3$, which due to (\ref{eq:proof_beta}) yields the condition (\ref{eq:theorem_qubic_condition}) of Theorem~\ref{theorem_cubic}. 
Substituting expressions for $\xi$, $\beta$ and the maximum value of $|b|$ from the condition (\ref{eq:theorem_qubic_condition}) in formulas (\ref{eq:proof_cubic_gamma}, (\ref{eq:proof_cubic_alpha}) for the parameters $\gamma$ and $\alpha$ we obtain boundaries (\ref{eq:theorem_qubic_domain}) of the invariant interval $D_x$ expressed via the parameter $\mu$ only. 

Therefore, similarly to the Theorem~\ref{theorem_logistic_Henon_attractor} the map $H$ with cubic nonlinearity $f(x)$ (\ref{eq:qubic}) has an attractor $A$, lying in the domain $D$ with boundaries $\alpha$ and $\gamma$ explicitly defined via the parameters of map~\eqref{map:T}, \eqref{eq:qubic}. $\Box$

\textbf{Remark.} With increase of the parameter $\gamma$ the auxiliary maps undergo the saddle-node bifurcation at the curve 
\begin{equation}
|b|=\dfrac{(1-a)(1+\mu)^{1.5}}{\mu^{1.5}+(1+\mu)^{1.5}}.
\label{eq:cubic_saddle-node}
\end{equation}
Comparing the inequality (\ref{eq:theorem_qubic_condition}) we conclude that the curve (\ref{eq:cubic_saddle-node}) is outside the parameters region (\ref{eq:theorem_qubic_condition}) for existence of the invariant domain $D$.

\section{Smale horseshoe of multidimensional H\'{e}non map}

Consider H\'{e}non map, i.e. the case~\eqref{map:T} for $f(x)=\mu-x^2$.
The auxiliary maps take the form
\begin{equation}
\bar{x}=\mu\pm\gamma-x^2\triangleq f^{\pm}(x).
\label{eq:aux_map_Henon}
\end{equation}

1) We choose interval $D_x=\lbrace |x|<\mu-\gamma\rbrace$
$(\alpha^{\pm}=\pm(\mu-\gamma))$.

2) Choose the region of parameters and intervals of $x$ when Theorem~\ref{th:attractor_gamma} does not work. 

That is
\begin{equation*}
     f(x)+\gamma>\alpha^+\quad\textrm{and}\quad
     f(x)-\gamma<\alpha^-,
\end{equation*}
which give
\begin{equation*}
          |x|<\sqrt{2\gamma}\quad\textrm{and}\quad 
          |x|>\sqrt{2(\mu+\gamma)},
\end{equation*}
respectively. 
These intervals lie at $D_x$ in the parameter region
\begin{equation}
    \sqrt{2(\mu+\gamma)}<\mu-\gamma.
    \label{eq:Smale_2sqrt<mu}
\end{equation}

\begin{theorem}
    In the parameter region \eqref{eq:Smale_2sqrt<mu} map~\eqref{map:T} has topological Smale horseshoe. 
    \label{th:Smale}
\end{theorem}

\textbf{Proof.} 
Consider arbitrary line $L=\lbrace y=\textrm{const}\,\in D_y, x\in D_x\rbrace$.
The image of $L$ is curve $(\bar{x}, \bar{y})=TL(x,y)$ given by formulas
\begin{equation*}
    \begin{array}{l}
         \bar{x}=f(x)+\textrm{const},\\
         \\
         \bar{y}_1=bx,\\
         \\
         \bar{y}_i=\textrm{const},\quad i=2,\ldots,n,
    \end{array}
\end{equation*}
where $x$ is parameter, $x\in(\alpha^-,\alpha^+)$.
This curve is an arc starting at $x=-(\mu-\gamma)$ with $\bar{x}(\alpha^-)<-(\mu-1)$, visiting region $\bar{x}>\mu-\gamma$ for $x\in(-\sqrt{2\gamma},\sqrt{2\gamma})$ and returning to region $\bar{x}<-(\mu-\gamma)$ for $x\in(\sqrt{2(\mu+\gamma)},\mu-\gamma)$.
Hence the whole set of lines $L$ filling in region $D_{\alpha}$ forms a solid arc which is due to construction is Smale horseshoe. $\Box$ 

\begin{remark}
    Using similar arguments one can prove the existence of Smale horseshoes for cubic or polynomial function $f(x)$.
\end{remark}

\section{The structure of attractors}

\subsection{Characteristic multipliers}

The Jacobi matrix $J$ of map~\eqref{map:T}, having the form
\begin{equation}
J=
\left(
\begin{array}{ccccccc}
f_x'&1&1&\cdots &1&1&1\\
b&0&0&\cdots &0&0&0\\
0&a_1&0&\cdots &0&0&0\\
\vdots&\vdots&\vdots&\ddots&\vdots&\vdots&\vdots\\
0&0&0&\cdots &a_{n-2}&0&0\\
0&0&0&\cdots &0&a_{n-1}&0
\end{array}
\right),
\label{eq:Jacobi_matrix}
\end{equation}
depends only on the $x$-coordinate. 
The eigenvalues of this matrix are the roots of the characteristic polynomial
\begin{equation}
\Delta_{n+1}=\det(J-sE)=0,
\end{equation}
which is the solution of the linear inhomogeneous difference equation
\begin{equation}
\Delta_{n+1}=(-1)^nb\prod\limits_{k=1}^{n-1}a_k-s\Delta_n,
\label{eq:Dn_diference_equation}
\end{equation}
where the determinant $\Delta_n$ takes the form
\begin{equation}
\Delta_n=
\left\vert
\begin{array}{ccccccc}
f_x'-s&1&1&\cdots &1&1&1\\
b&-s&0&\cdots &0&0&0\\
0&a_1&-s&\cdots &0&0&0\\
\vdots&\vdots&\vdots&\ddots&\vdots&\vdots&\vdots\\
0&0&0&\cdots &a_{n-3}&-s&0\\
0&0&0&\cdots &0&a_{n-2}&-s
\end{array}
\right\vert .
\label{eq:Dx_matrix}
\end{equation}

The initial condition for this equation is the determinant 
\begin{equation}
\Delta_2=s^2-f_x's-b, 
\end{equation}
which serves the characteristic polynomial to the original 2-D H\'{e}non map. 
With this initial condition the solution of the equation (\ref{eq:Dn_diference_equation}),(\ref{eq:Dx_matrix}) has the form (see Appendix A)
\begin{equation}
(-1)^{n+1}\Delta_{n+1}=s^{n-1}(s^2-f_x's-b)-
b\sum\limits_{k=3}^{n+1}\left(\prod\limits_{l=1}^{k-2}a_l\right)s^{n-k+1}.
\label{eq:D_n+1_solution}
\end{equation}

Note, that the sum in the solution (\ref{eq:D_n+1_solution}) has the simple form
\begin{equation}
\sum=a_1s^{n-2}+a_1a_2s^{n-3}+\cdots+a_1a_2\ldots a_{n-1}.
\end{equation}

\subsection{The case of small parameter b}

Consider the case of small parameter $|b|$. 
The next statement holds.

\begin{theorem}
Let the parameter $|b|$ of map~\eqref{map:T} be small. 
Then each structurally stable $p$-periodic orbit of the 1-D map $x\rightarrow f(x)$ defines the similar $p$-periodic orbit of map~\eqref{map:T}. 
\label{theorem_b_small}
\end{theorem}

\textbf{Proof.} 
In the limiting case $b=0$ map~\eqref{map:T} reduces to the form 
\begin{equation}
\begin{array}{lcl}
\bar{x}=&f(x)+\mathds{1}y&\triangleq \; g(x,y),\\
\bar{y}=&Ay&\triangleq \; L(0,y),
\end{array}
\label{sys:H_b=0}
\end{equation}

The map $L(0,y)$, being independent on $x$, due to obvious inequality $\Vert \bar{y} \Vert <a\Vert y\Vert$
(see (\ref{eq:a<1_b<1}), (\ref{eq:norm_ineq})) is contractive and has the unique absolutely stable fixed point $y=0$.
In fact, the dynamics of the map $L(0,y)$ is such that for any initial point $y_0$ at the first iteration the coordinate $y_1$ becomes zero, at the second -- $y_2$ becomes zero, and so that in $(n-1)$ iterations the vector $y$ reaches the zero fixed point. 
Such behaviour of $y$ coordinates is due to $n$ zero eigenvalues of the matrix $A$.
Hence, the dynamics of the map (\ref{sys:H_b=0}) is defined by 1-D map $g(x,0): x\rightarrow f(x)$.

One can obtain the same statement using the characteristic polynomial (\ref{eq:D_n+1_solution}) for $b=0$
\begin{equation}
\Delta_{n+1}=s^n(s-f_x').
\label{eq:proof_theorem_b_small_char_polyn}
\end{equation}

Let the map $f(x)$ has a $p$-periodic orbit $\lbrace x_1, x_2,\ldots,x_p\rbrace$, which corresponds to the periodic orbit $O_0=\lbrace[x(1),0],[x(2),0],\ldots, [x(p),0]\rbrace$ of the map (\ref{sys:H_b=0}).
The orbit $O_0$  has $n$ zero multipliers corresponding to $y$ variables and the multiplier for $x$ variable
\begin{equation}
m_x=\prod\limits_{i=1}^p f_x'(x(i)).
\label{eq:m_x_multiplier}
\end{equation}

We assume that $m_x\neq 0$, $|m_x|\neq 1$.
As far as multipliers do not lye on the unit circle the periodic orbit $O_0$ is structurally stable.

Note, that the reduced 1-D map $f(x)=g(x,0)$ in general case is non-invertible,
but in a small vicinity of the periodic orbit $O_0$ this map is invertible with respect to certain kneading corresponding to this orbit. 
Using this property we obtain that when $|b|$ increases from zero this orbit persists due to its structure stability. 
This orbit leaves the manifold $y=0$ and becomes the orbit $O_y=[(x(1),y(1)),(x(2),y(2)),\ldots, (x(p),y(p))]$ of the invertible map~\eqref{map:T}.
The orbit $O_y$ has the multipliers which are close to those of the orbit $O_0$ of the map (\ref{sys:Henon_1}) due to continuous and smooth dependence on the parameters of the roots of the polynomial $(\ref{eq:D_n+1_solution})$. $\Box$

\section{Acknowledgements}
This work was supported by the Russian Science Foundation under grant No. 22-21-00553.

\bibliographystyle{unsrt}

\end{document}